\documentclass{amsart}

\usepackage[english]{babel}

\usepackage[letterpaper,top=2cm,bottom=2cm,left=3cm,right=3cm,marginparwidth=1.75cm]{geometry}

\newtheorem{theorem}{Theorem}[section]
\newtheorem{lemma}[theorem]{Lemma}

\newtheorem{problem}[theorem]{Problem}

\theoremstyle{definition}
\newtheorem{definition}[theorem]{Definition}
\newtheorem{example}[theorem]{Example}

\theoremstyle{remark}

\usepackage{amsmath,amssymb}
\usepackage{graphicx}
\usepackage[colorlinks=true, allcolors=blue]{hyperref}
\usepackage{todonotes}
\usepackage{tikz}
\usetikzlibrary{decorations.pathmorphing}
\usepackage{ytableau}
\usepackage{verbatim}
\usepackage{diagbox}

\newcommand{\cF}{\mathcal{F}}

\newcommand{\B}{\mathcal{B}}
\renewcommand{\L}{\mathcal{L}}

\newcommand{\rank}{\operatorname{rank}}

\newcommand{\RR}{\mathbb{R}}

\newcommand{\Le}{\mathsf{Le}}
\newcommand{\Sq}{\mathsf{Sq}}

\newcommand{\sm}{\setminus}

\DeclareMathOperator{\Meas}{Meas}
\DeclareMathOperator{\supp}{supp}
\DeclareMathOperator{\Gr}{Gr}

\DeclareMathOperator{\Net}{Net}

\newcommand{\Z}{\mathbb{Z}}
\usepackage{xspace}
\newcommand{\bv}{\mathbf{v}}
\newcommand{\be}{\mathbf{e}}

\title{Transversal Positroids}
\author{John Machacek \and George D. Nasr}
\address{Department of Mathematics and Computer Science, Hampden-Sydney College} \email{jmachacek@hsc.edu}
\address{Department of Mathematics, Augustana University} \email{george.nasr@augie.edu}

\begin{document}
\maketitle

\begin{abstract}
We exhibit three infinite families of positroids that are transversal matroids.
The first comes from a subclass of Postnikov's $\Le$-diagrams that we call $\Sq$-diagrams.
Our other two families are obtained by classifying all sparse paving positroids along with all rank $2$ positroids that are transversal.
Subsequently, we construct transversal sparse paving positroids that lack a noncrossing minimal presentation, disproving a conjecture of Marcott.
\end{abstract}

\section{Introduction}
Positroids and transversal matroids are two types of real representable matroids.
This means both can be represented by matrices with real entries.
Our focus is on the study of the overlap between these classes of matroids.
Positroids come from real-valued matrices for which all maximal minors are nonnegative.
We will often work with a particular matrix coming from weighted paths in a planar network.
Meanwhile, transversal matroids come from real-valued matrices with a fixed pattern of zero entries along with ``generic'' nonzero entries.
Alternatively, transversal matroids can be presented by systems of sets or equivalently as matchings in bipartite graphs.
We will formally define positroids and transversal matroids in Section~\ref{sec:prelim}, along with other necessary background for this paper.

Positroids play a role in certain computations in $\mathcal{N}=4$ SYM theory in high-energy physics (see e.g.,~\cite{Amp}).
Transversal positroids in particular appear as the matroids associated to Wilson loop diagrams~\cite{WilsonLoop}.
Also, transversal positroids arise in the study of the amplituhedron~\cite{KWZ} (see Example~\ref{ex:amp}).
These connections to physics were one motivation for Marcott's systematic study of transversal positroids~\cite{MarcottThesis} where the focus was on Problem~\ref{prob:TtoP}.

\begin{problem}
Determine if a given transversal matroid is a positroid.
\label{prob:TtoP}
\end{problem}

Marcott has made progress on Problem~\ref{prob:TtoP} and conjectured a complete solution in terms of a noncrossing condition~\cite[Theorem 3.4.4 and Conjecture 3.4.5]{MarcottThesis}.
In Theorem~\ref{thm:Sq}, for $\Sq$-diagrams defined in Definition~\ref{def:Sq}, we connect Marcott's noncrossing condition to the topological noncrossing of planar networks used by Postnikov~\cite{Pos} to parameterize the totally nonnegative Grassmannian.
However, Theorem~\ref{thm:counter} provides counterexamples to Marcott's conjecture coming from certain sparse paving positroids.
Most of our work relates to the reverse of Problem~\ref{prob:TtoP} where one tries to decide if a positroid is transversal.
We are also interested in figuring out when a positroid is a special type of transversal matroid known as a fundamental transversal matroid.

\begin{problem}
Determine if a given positroid is a (fundamental) transversal matroid.
\label{prob:PtoT}
\end{problem}

In general, deciding if a matroid is transversal can be a tricky problem.
In~\cite[Theorem 1 (11)]{JensenKorte}, it is shown that there is no polynomial-time algorithm to decide whether a matroid is transversal using an independence test oracle (see also~\cite[Example 3.12]{VW}).
Also, even given a transversal matroid there is no polynomial-time algorithm to decide whether it is fundamental transversal using an independence test oracle~\cite[Theorem 1 (12)]{JensenKorte}.
However, we have characterizations of transversal and fundamental transversal matroids involving inequalities of ranks~\cite{BoninKungdeMier}.

In approaching Problem~\ref{prob:PtoT}, we must choose how to represent positroids.
We will primarily work with the $\Le$-diagrams and $\Le$-Networks of Postnikov~\cite{Pos}, and are able to solve this problem for some infinite families.
In one case, we do this by defining a subclass of $\Le$-diagrams that we call $\Sq$-diagrams in Definition~\ref{def:Sq} and show that they always produce transversal positroids in Theorem~\ref{thm:Sq}.
The case of sparse paving matroids that are also positroids is solved in Theorem~\ref{thm:sparsepaving} using a characterization of sparse paving positroids from~\cite{MachacekNasr}.
In Theorem~\ref{thm:rk2} we give a simple criterion on the $\Le$-diagram of a rank $2$ positroid that determines when the positroid is a transversal matroid.

\section{Preliminaries}
\label{sec:prelim}

\subsection{Matroids}
Let $E$ be a finite set and let $2^E$ denote the set of all subsets of $E$.
A \emph{matroid basis system on ground set $E$} is a nonempty family $\B \subseteq 2^E$ of \emph{bases} satisfying the \emph{exchange axiom}: for all distinct $B,B'\in\B$ and all $e\in B\sm B'$, there exists $e'\in B'\sm B$ such that $\left(B\sm \{e\}\right)\cup\{e'\}\in \B$.
A \emph{matroid} is defined by the pair $M=(E,\B)$.
Any $I \subseteq B \in \B$ is an \emph{independent set}.
A matroid's \emph{rank function} is given by $\rank(A)=\max\{|A\cap B|:\ B\in\B\}$ for any $A\subseteq E$.
The \emph{rank} of $M$ can be denoted $\rank(M)$ and is equal to $\rank(E)$.
A subset of $E$ that is not independent is called \emph{dependent}.
A \emph{circuit} is a minimal, by inclusion, dependent set.
A subset $F\subseteq E$ is called \emph{flat} of $M$ when $\rank(F) < \rank(F \cup \{x\})$ for any $x \in E \sm F$.
When $F$ is a flat of $M$ with $\rank(F) = \rank(M)-1$ we call $F$ a \emph{hyperplane}.
A \emph{cyclic flat}  of a matroid $M$ is a flat of M that is the union of a (possibly empty) set of circuits.
There always exists the maximal cyclic flat (the set of non-coloops) and the minimal cyclic flat (the set of loops).

Given a matroid $M$ we have its dual matroid $M^*$ on the same ground set $E$ whose bases are complements in $E$ of bases of $M$.
A \emph{loop} is a circuit that is a singleton.
Dual to a loop is the notion of a \emph{coloop}, which is an element $e$ that is contained in every basis.
\emph{Cocircuits} are circuits in the dual matroid, or equivalently, they are complements of hyperplanes.

In general, proving that a collection of sets is a matroid basis system is not a trivial task, though some examples require less work than others.
For example, take a nonnegative integer $k$ and finite set $E$, the \emph{uniform matroid} $U_{k}(E)$ on ground set $E$ has basis system $\binom{E}{k}:=\{B\subseteq E:\ |B|=k\}$. 
We leave it to the reader to verify that this collection of sets indeed satisfies the exchange axiom.
We will often find it convenient to write $U_{k,n}:=U_k([n])$, where $[n]:=\{1,2,\dots, n\}$.

Given a full rank\footnote{It is possible to use any field, and the choice of the field can affect the matroid. In this work, we will always use the field of real numbers since we are concerned with positivity.} $k \times n$ matrix $A$ one can obtain a \emph{linear matroid} denoted $M_A := ([n], \B)$ where $\B \subseteq \binom{[n]}{k}$ consists of those $k$-element subsets of columns that index nonvanishing maximal minors of $A$.
Equivalently, $\B$ consists of collections of column vectors from $A$ that form a basis of $\RR^k$.
Linear matroids are also called \emph{representable matroids}.

\subsection{Transversal Matroids}\label{sec:trans}
A \emph{set system} $(S_1, S_2, \dots, S_k)$ is a family of subsets of some ground set $E$.
A set system may contain repeated sets; that is, we may have $S_i = S_j$ for $i \neq j$.
Also, even though they are written as a tuple, we do not distinguish between rearrangements of a set system.
A \emph{transversal} of a set system $(S_1, S_2, \dots, S_k)$ is a set $I\subseteq E$ for which there is a bijection $\phi:[k] \to I$ so that $\phi(j) \in S_{j}$ for all $j \in [k]$.
We will assume that some transversal exists for all set systems we consider.
This assumption is consistent with the motivating work of Marcott~\cite{MarcottThesis}.
The \emph{transversal matroid} $M[(S_1, S_2, \dots, S_k)]$ is the matroid whose bases are the transversals of the set system $(S_1, S_2, \dots, S_k)$.
So, we say that $(S_1, S_2, \dots, S_k)$ is a \emph{presentation} of the transversal matroid $M = M[(S_1, S_2, \dots, S_k)]$.
Our assumption that a transversal exists means we will only consider set systems with $k$ sets when presenting a rank $k$ transversal matroid.
Though it is possible to consider more general presentations, it is always possible to find a presentation of a rank $k$ transversal matroid with a set system with $k$ sets~\cite[Lemma 2.4.1]{Oxley}.
Thus, our assumption on presentations does not limit which transversal matroids we consider.

In general, many set systems give us the same transversal matroid.
We may define a partial order denoted $\preceq$ on presentations of a given transversal matroid by saying $(S'_1, S'_2, \dots, S'_k) \preceq (S_1, S_2, \dots, S_k)$ when (for some rearrangement) $S'_i \subseteq S_i$ for each $i$ and $M[(S'_1, S'_2, \dots, S'_k)] = M[(S_1, S_2, \dots, S_k)]$.
A \emph{minimal presentation} is a presentation that is minimal with respect to $\preceq$.
A transversal matroid can have many minimal presentations.

\begin{example}\label{ex:uniform}
    Every uniform matroid $U_{k,n}$ is a transversal matroid. For instance, we have
    \[U_{k,n}=M[(\underbrace{[n],[n],\dots, [n]}_{k})].\]
    Another presentation for this matroid is
    \[U_{k,n}=M[([1,n-k+1],[2,n-k+2],\dots, [k,n])],\]
    where for integers $i$ and $j$, we have $[i,j]:=\{k\in \Z: i\leq k \leq j\}$ denoting an \emph{interval} in $\Z$. 
\end{example}

All transversal matroids are linear matroids, and they can be represented using a matrix of real numbers~\cite[Corollary 11.2.7]{Oxley}.
The correspondence goes as follows.
Consider $(S_1, S_2, \dots, S_k)$ where $S_i \subseteq [n]$ for each $i$.
Then build a $k \times n$ matrix $A$ where the $i$th row is a vector $v$ so that $v_j = 0$ whenever $j \not\in S_i$ and $v_j$ is a generic\footnote{If one prefers, the nonzero entries can be taken to be distinct formal variables which are later specialized to real numbers} nonzero entry whenever $j \in S_i$.

\begin{example}
If we take the set system
\[(\{1,4,5\}, \{2,4,5\}, \{3,5\})\]
we would get the matrix
\[\begin{bmatrix}
    * & 0 & 0 & * & *\\
    0 & * & 0 & * & *\\
    0 & 0 & * & 0 & *
\end{bmatrix}\]
where an $*$ denotes a generic real entry.
One can check that the transversals of the set system are exactly the collections of three columns that form a nonsingular $3\times 3$ matrix.
For example, $\{1,4,5\}$ is a transversal, and
\[\det\left(\begin{bmatrix} * & * & * \\ 0 & * & *\\ 0 & 0 & *\end{bmatrix}\right) \neq 0\]
for any generic choice of the nonzero entries.
Also, $\{1,2,4\}$ is not a transversal, and
\[\det\left(\begin{bmatrix} * & 0 & * \\ 0 & * & *\\ 0 & 0 & 0\end{bmatrix}\right) = 0\]
for any choice of the nonzero entries.
\end{example}

A \emph{fundamental transversal matroid} is a transversal matroid $M[(S_1, \dots, S_k)]$ on ground set $E$ such that for each $1 \leq i \leq k$ there exists $e \in E$ where $e \in S_i$ and $e \not\in S_j$ for $j \neq i$.
This is one of many descriptions of fundamental matroids (see e.g.,~\cite{BoninKungdeMier}).

\subsection{The Grassmannian and positroids}
\label{sec:Grass}
The \emph{(real) Grassmannian} $\Gr(k,n)$ is the set of $k$-dimensional subspaces of $\RR^n$.
We may assume $0 \leq k \leq n$.
Any element of $V \in \Gr(k,n)$ can be represented by a full-rank $k \times n$ matrix $A$ so that $V$ is the row span of $A$.
In this case we write $V = [A]$ and note that the choice of $A$ is unique up to left multiplication by elements of the \emph{general linear group} $GL_k$ consisting of invertible $k \times k$ real matrices.

Recall, $[n] = \{1,2,\dots, n\}$ and so $\binom{[n]}{k}$ denotes the collection of all $k$-element subsets of $[n]$.
Given $I \in \binom{[n]}{k}$ and a $k \times n$ matrix $A$, we let $\Delta_I(A)$ denote the value of the maximal minor with columns indexed by $I$.
The \emph{totally nonnegative Grassmannian} is 
\[\Gr^{\geq 0}(k,n) = \left\{V \in \Gr(k,n) : V = [A] \text{ and } \Delta_I(A) \geq 0 \text{ for all } I \in \binom{[n]}{k}\right\}\]
as defined by Postnikov~\cite{Pos}.
A \emph{positroid} is then a linear matroid $M_A$ such that $[A] \in \Gr^{\geq 0}(k,n)$.

We always consider positroids with an ordered ground set.
This convention is typical when considering the Grassmannian and total positivity.
The ordering corresponds to the order of columns in a matrix representing a linear subspace, and reordering the columns can change the sign of the maximal minors.
With this approach, it is possible that positroids that differ from our perspective can be isomorphic as abstract matroids.
Of course, determining whether a positroid is transversal depends only on its isomorphism class.
However, with an ordered ground set, distinct positroids correspond to distinct cells in Postnikov's~\cite {Pos} cell decomposition of the totally nonnegative Grassmannian, and we wish to be able to consider any cell in this decomposition.
We will now explain a construction due to Postnikov from which we can obtain matrices representing all possible positroids.

\begin{definition}[Network]
A \emph{planar directed network}, or simply \emph{network}, is a planar directed acyclic\footnote{The acyclic assumption is not necessary, but it is sufficient for our purpose.} graph $N = (V, E)$ with embedding in a disk so that:
\begin{itemize}
    \item $[n] \subseteq V$ and the elements of $[n]$ are \emph{boundary vertices}  embedded on the boundary of the disk in order, going clockwise.
    \item each boundary vertex is designated a source or sink and is incident on at most one edge.
    \item each directed edge $e$ comes with a strictly positive real edge weight $w_e$.
\end{itemize}
We let $\Net(k,n)$ denote the collection of planar directed networks with $n$ boundary vertices of which $k$ are sources.
\end{definition}

We allow for isolated sources or sinks.
It is important to know whether an isolated boundary vertex is a source or a sink, as it will correspond to a coloop or loop, respectively.
An example network is shown on the right in Figure~\ref{fig:LeNetwork}.
Edges that do not have a weight shown are, by convention, assumed to have a weight equal to $1$.
We indicate that the boundary vertex $7$ is a sink by including an edge directed towards it that does not originate at any vertex in the network.

\begin{definition}[Boundary Measurement]
Given a network $N$ with $K \subseteq [n]$ being the set of boundary sources, we get the $K \times [n]$ \emph{boundary measurement matrix} of $N$ defined by
\[\Meas(N)_{ij} = \sum_{P: i \rightsquigarrow j}(-1)^{s_{ij}} \prod_{e \in P} w_e\]
for each $i \in K$ and $j \in [n]$.
Here the $P: i \rightsquigarrow j$ denotes a directed path from $i$ to $j$, and $e \in P$ is a directed edge along this path.
Note the acyclic assumption ensures each entry of the matrix is a finite sum.
The value $s_{ij}$ counts the number of sources strictly between $i$ and $j$ in the usual order.
\end{definition}

For $i \in K$ we have $\Meas(N)_{ii} = 1$, which can be thought of as coming from the empty path $i \rightsquigarrow i$.
This construction is due to Postnikov, who proved that $\Meas(N)$ represents a positroid and every positroid can be represented by some network~\cite[Theorem 4.8]{Pos}.
There is another object related to positroids we will need, which we now define.

\begin{definition}[{$\Le$-diagram~\cite[Definition 6.1]{Pos}}]
For any partition $\lambda$ we consider its Young diagram with the convention that our diagrams are left-justified, with the largest part on top.
We index boxes of Young diagrams by pairs $(i,j)$ such that $i$ increases as we go downward and $j$ increases as we go right.
A \emph{$\Le$-diagram} $D$ is a filling of the Young diagram so that each box is either empty or contains $\bullet$ such that for $i < i'$ and $j < j'$  the box $(i',j')$ exists it contains $\bullet$ whenever each the boxes $(i,j')$ and $(i',j)$ contain a $\bullet$.
We let $|D|$ denote the number of $\bullet$'s that $D$ contains, and we let $\Le(k,n)$ denote the collection of all $\Le$-diagrams that fit inside a $k \times (n-k)$ box.
\end{definition}

A key feature of $\Le$-diagrams is that they can readily be turned into networks~\cite[Definition 6.3 and Theorem 6.5]{Pos}.
Given $D \in \Le(k,n)$, the southeast boundary of the shape, along with possibly the top and left sides of the box, gives us a lattice path from the top right to the bottom left of the box.
We will start in the upper right corner of the $k \times (n-k)$ box and place a vertex along each step, labeled with an element from $[n]$ in order.
Vertical steps are sources, and horizontal steps are sinks.
Any $\bullet$ in $D$ is an internal vertex.
A directed edge is given a positive real weight and drawn from right to left from each source or internal vertex whenever there is an internal vertex to the left in the same row.
Also, a directed edge is drawn down from an internal vertex whenever there is an internal vertex or sink below in the same column.
A weight equal to $1$ is given to all downward vertical edges.

\begin{figure}
\begin{minipage}{0.48\textwidth}

\centering
\ytableausetup{centertableaux}
\begin{ytableau}
$ $&\bullet& \bullet & $ $ \\
$ $&\bullet & \bullet & \bullet\\
$ $& & \bullet \\
\end{ytableau}
\end{minipage}\hfill
\begin{minipage}{0.48\textwidth}
\begin{tikzpicture}
        \draw  (0,0) circle (2cm) {};
         \node[circle,fill=black, scale=0.4, label=right:{$1$}] (1) at ({2*cos(50)}, {2*sin(50)}) {};
         \node[circle,fill=black, scale=0.4, label=right:{$2$}] (2) at ({2*cos(20)}, {2*sin(20)}) {};
         \node[circle,fill=black, scale=0.4, label=right:{$3$}] (3) at ({2*cos(-10)}, {2*sin(-10)}) {};
         \node[circle,fill=black, scale=0.4, label=right:{$4$}] (4) at ({2*cos(-35)}, {2*sin(-35)}) {};
         \node[circle,fill=black, scale=0.4, label=below:{$5$}] (5) at ({2*cos(-90)}, {2*sin(-90)}) {};
         \node[circle,fill=black, scale=0.4, label=below:{$6$}] (6) at ({2*cos(-120)}, {2*sin(-120)}) {};
        \node[circle,fill=black, scale=0.4, label=below:{$7$}] (7) at ({2*cos(-160)}, {2*sin(-160)}) {};

        \node[circle,fill=black, scale=0.4] (21) at ({cos(20)}, {2*sin(20)}) {};
        \draw[-{latex}] (2) -- (21) node[midway,below] {$c$};
        \draw[-{latex}] (21) to[bend right] (3);

        \node[circle,fill=black, scale=0.4] (12) at ({2*cos(-90)}, {2*sin(50)}) {};
         \draw[-{latex}] (1) -- (12)  node[midway,above] {$a$};
         \node[circle,fill=black, scale=0.4] (22) at ({2*cos(-90)}, {2*sin(20)}) {};
         \draw[-{latex}] (21) -- (22) node[midway,below] {$d$};
         \node[circle,fill=black, scale=0.4] (32) at ({2*cos(-90)}, {2*sin(-35)}) {};
         \draw[-{latex}] (4) -- (32) node[midway,below] {$f$};
         \draw[-{latex}] (12) -- (22);
        \draw[-{latex}] (22) -- (32);
        \draw[-{latex}] (32) -- (5);

    \node[circle,fill=black, scale=0.4] (23) at ({2*cos(-120)}, {2*sin(20)}) {};
    \node[circle,fill=black, scale=0.4] (13) at ({2*cos(-120)}, {2*sin(50)}) {};
    \draw[-{latex}] (12) -- (13) node[midway,above] {$b$};
    \draw[-{latex}] (13) -- (23);
    \draw[-{latex}] (22) -- (23) node[midway,below] {$e$};
    \draw[-{latex}] (23) -- (6);

 \end{tikzpicture}
 \end{minipage}
\caption{A $\Le$-diagram $D$ (left) and its network $N_D$ (right).}
\label{fig:LeNetwork}
\end{figure}

\begin{example}[$\Le$-diagram, network, and boundary measurement]
In Figure~\ref{fig:LeNetwork} we find both a $\Le$-diagram $D$ on the left and its network $N_D$ on the right.
The boundary vertex $7$ plays the role of a sink in this network since it comes from the horizontal step in the border of the $\Le$-diagram. The boundary measurement matrix is then
\[\Meas(N_D) = 
\begin{bmatrix}
1 & 0 & 0 & 0 & a & ab+ae & 0 \\
0 & 1 & c & 0 & -cd & -ced & 0\\
0 & 0 & 0 & 1 & f & 0 & 0
\end{bmatrix}\]
where the rows are indexed by the boundary vertices $1$, $2$, and $4$.
\label{ex:ND}
\end{example}

\subsection{Noncrossing presentations}
We will be interested in a special subclass of presentations of transversal matroids, which we define now. For each $a \in [n]$ we define the order $<_a$ by $a <_a a+1 \cdots <_a n <_a 1 <_a 2 \cdots <_a a-1$, which is simply a cyclic shift of the usual total order on $[n]$.
For $S,T \subseteq [n]$, we say $S$ \emph{crosses} a set $T$ if there exists $a,b,c,d \in [n]$ such that:
\begin{itemize}
    \item[(i)] $a <_a b <_a c <_a d$,
    \item[(ii)] $a,c \in S$ and $a,c \not\in T$,
    \item[(iii)] $b,d \in T$ and $b \not\in S$.
\end{itemize}
This definition is not symmetric.
In other words, $S$ crossing $T$ is not an equivalent condition to $T$ crossing $S$.
A set system $(S_1, S_2, \dots, S_k)$ is called \emph{crossing} if there exisits $i$ and $j$ so that $S_i$ crosses $S_j$, otherwise the set system is called \emph{noncrossing}.
Marcott~\cite[Theorem 3.4.4]{MarcottThesis} has shown that if $(S_1, S_2, \dots, S_k)$ is a noncrossing minimal presentation, then $M[(S_1, S_2, \dots, S_k)]$ is positroid.
Furthermore, Marcott~\cite[Conjecture 3.4.5]{MarcottThesis} has conjectured that there exists a noncrossing minimal presentation for any positroid that is also a transversal matroid.

In some instances, it will be convenient to think about our set systems as consisting of complements.
From this perspective we find that $[n]\sm X$ crosses $[n] \sm Y$ when there exists $a,b,c,d \in [n]$ with
\begin{itemize}
    \item[(i)] $a <_a b <_a c <_a d$,
    \item[(ii)] $a,c \not\in X$ and $a,c \in Y$,
    \item[(iii)] $b,d \not\in Y$ and $b \in X$.
\end{itemize}
which follows readily from the conditions above.

For any $0 \leq k \leq n$ and $i \in [n]$, we define $C_{k,n}^{(i)} \in \binom{[n]}{k}$ as $C_{k,n}^{(i)} = \{i, i+1, \dots i + k -1\}$ where we work modulo $n$ using $[n]$ as the set of representatives.
Hence, $C_{k,n}^{(i)}$ is the \emph{cyclic\footnote{This cyclic comes from wrapping around modulo $n$ and is unrelated to matroid terminology of cyclic flats.} interval} of size $k$ in $[n]$ that starts at $i$.
Cyclic intervals will be of particular importance to us, and we now make a few observations about them.
First, the complement of a cyclic interval in $[n]$ is another cyclic interval.
Secondly, if $S$ crosses $T$, then $T$ cannot be a cyclic interval.
This is because having $a <_a b <_a c <_a d$ with $a,c \not\in T$ and $b,d \in T$ contradicts the structure of a cyclic interval.

\section{A class of fundamental transversal positroids}

Given any vector $v = (v_1, \dots, v_n) \in \RR^n$ we have its \emph{support} given by
\[\supp(v) := \{i : v_i \neq 0 \text{ and }1 \leq i \leq n\}\]
which simply consists of the indices of nonzero entries of the vector.
For a $k \times n$ matrix $A$ with rows $r_i \in \RR^n$ for $1 \leq i \leq k$ we define
\[\supp(A) := (\supp(r_1), \supp(r_2), \dots, \supp(r_k))\]
which we may view as a set system.

\begin{lemma}
If $N \in \Net(k,n)$, then $\supp\left(\Meas(N)\right)$ is noncrossing.
\label{lem:MeasNoncrossing}
\end{lemma}

\begin{figure}[h]
    \centering
    \begin{tikzpicture}
    %%% LEFT NETWORK
    \draw (0,0) circle (2.5cm);
    \node (a) at (0,2.5) {};
    \node (aa) at (0,2.7) {$a$};
    
    \node (b) at (2.5,0) {};
    \node (bb) at (2.7,0) {$b$};
     
    \node (c) at (0,-2.5) {};
    \node (cc) at (0,-2.7) {$c$};
     
    \node (d) at (-2.5,0) {};
    \node (dd) at (-2.7,0) {$d$};
    
    \node (S) at ({2.5*cos(225)}, {2.5*sin(225)}) {};
    \node (SS) at ({2.5*cos(225)-0.2}, {2.5*sin(225)-0.2}) {$S_i$};
    
    \node[circle, fill=gray!, scale=0.4] (O) at (0,0) {};
    \draw[decorate, decoration={snake},-{latex}] ({2.5*cos(225)}, {2.5*sin(225)}) -- (O);
    \draw[decorate, decoration={snake},-{latex}]  (O) -- (0,2.5);
    \draw[decorate, decoration={snake},-{latex}]  (O) -- (0,-2.5);

     %%% RIGHT NETWORK
    \begin{scope}[shift={(7,0)}]
    \draw (0,0) circle (2.5cm);
    \node (a) at (0,2.5) {};
    \node (aa) at (0,2.7) {$a=S_i$};
    
    \node (b) at (2.5,0) {};
    \node (bb) at (2.7,0) {$b$};
     
    \node (c) at (0,-2.5) {};
    \node (cc) at (0,-2.7) {$c$};
     
    \node (d) at (-2.5,0) {};
    \node (dd) at (-2.7,0) {$d$};
    
    \draw[decorate, decoration={snake},-{latex}]  (0,2.5) -- (0,-2.5);

    \end{scope}

    \end{tikzpicture}
    \caption{Situations used in proof of Lemma~\ref{lem:MeasNoncrossing} that prevent crossings in set systems coming from supports of boundary measurement matrices.}
    \label{fig:NoCrossing}
\end{figure}

\begin{proof}
Take $N \in \Net(k,n)$ and let $A = \Meas(N)$  with rows $r_i \in \RR^n$ for $1 \leq i \leq k$.
Assume that $\supp(A)$ has a crossing so that $\supp(r_i)$ crosses $\supp(r_j)$ for some $i$ and $j$.
This means there exists $a$, $b$, $c$ and $d$ so that
\begin{align*}
  a,c \in \supp(r_i) &&  b,d \in \supp(r_j)\\
  b \not\in \supp(r_i) && a,c, \not\in \supp(r_j) 
\end{align*}
where $a <_a b <_a < c <_a d$.

Let $S_i$ and $S_j$ be the sources indexing the rows $r_i$ and $r_j$ respectively.
This means we must have paths $S_i \rightsquigarrow a$ and $S_i \rightsquigarrow c$.
It could be the case that both $a$ and $c$ are sinks which is shown on the left in Figure~\ref{fig:NoCrossing}.
Here the placement of $S_i$ can vary, but the fact that there is a continuous (undirected) curve from $a$ to $c$ will not change.
Indeed, we may always choose a path $S_i \rightsquigarrow a$ and a path $S_i \rightsquigarrow c$ which begin overlapping at $S_i$ then eventually diverge (if they intersect after diverging we can simply follow an overlapping portion for longer).
Alternatively, one of $a$ or $c$ could be a source.
The situation where $a$ is a source is shown on the right of Figure~\ref{fig:NoCrossing}.

We find that, in any case, we cannot have $S_j \rightsquigarrow b$ and $S_j \rightsquigarrow d$ without having a path departing from $S_j$ intersect the (undirected) continuous curve between $a$ and $c$.
Such an intersection would force a path $S_j \rightsquigarrow a$ or $S_j \rightsquigarrow c$.
This is a contradiction to our assumption that $\supp(r_i)$ crosses $\supp(r_j)$.
Therefore $\supp(A)$ must be noncrossing.
\end{proof}

\begin{figure}
\centering
\ytableausetup{centertableaux}
\begin{ytableau}
$ $ & & & & & & & & & &  & & &\\
$ $ & $d$ & & & & $a$ & & & & \\
& & & & &  & & &\\
& & &  & & &\\
$ $ & $c$ & & & & $b$ \\
 & & & &\\
& &
\end{ytableau}
    \caption{When boxes labeled $a$ and $c$ are each nonempty, then to be a $\Le$-diagram, the box labeled $b$ must be nonempty provided it exists. When boxes labeled $a$, $b$, and $c$ are all nonempty, then to be a $\Sq$-diagram, the box labeled $d$ also must be nonempty.}
    \label{fig:LeSq}
\end{figure}

\begin{definition}[$\Sq$-diagram].
A \emph{$\Sq$-diagram} $D$ is a $\Le$-diagram with the additional requirement that for $i < i'$ and $j < j'$ the box $(i,j)$ contains $\bullet$ whenever each the boxes $(i,j')$, $(i',j)$, and $(i',j')$ contain $\bullet$.
Note that the condition that $(i',j')$ contains a $\bullet$ implies this box must exist.
We let $\Sq(k,n)$ denote the collection of all $\Sq$-diagrams that fit inside a $k \times (n-k)$ box.
\label{def:Sq}
\end{definition}

The collection of $\Sq$-diagrams is smaller than the collection of all $\Le$-diagrams.
However, there are infinite families of interest consisting completely of $\Sq$-diagrams.
The following example shows the occurrence of certain $\Sq$-diagrams in the study of the amplituhedron.

\begin{example}[$m=2$ amplituhedron cells]
A collection of $\Sq$-diagrams that index a triangulation of the $m=2$ amplituhedron was given in~\cite[Definition 4.6]{KWZ}.
For given $k$ and $n$, this collection contains $\Sq$-diagrams of shape $\lambda = (\lambda_1, \dots, \lambda_k)$ that fit inside a $k \times (n-k)$ rectangle with $\lambda_i \geq 2$ for $1 \leq i \leq k$.
Boxes $(i,1)$ and $(i,\lambda_i)$ contain a $\bullet$ for $1 \leq i \leq k$ and all other boxes are empty.
The diagrams for $n=6$ and $k=3$ are shown in Figure~\ref{fig:m=2}.
\label{ex:amp}
\end{example}

\begin{figure}[h]
\centering
\ytableausetup{centertableaux}
\begin{ytableau}
\bullet &  & \bullet \\
\bullet &  & \bullet \\
\bullet & & \bullet
\end{ytableau}
\hspace{0.5cm}
\begin{ytableau}
\bullet &  & \bullet \\
\bullet &  & \bullet \\
\bullet & \bullet
\end{ytableau}
\hspace{0.5cm}
\begin{ytableau}
\bullet &  & \bullet \\
\bullet &  \bullet \\
\bullet & \bullet
\end{ytableau}
\hspace{0.5cm}
\begin{ytableau}
\bullet &   \bullet \\
\bullet &   \bullet \\
\bullet &  \bullet
\end{ytableau}
    \caption{Diagrams for top dimensional cells in $m=2$ amplituhedron for $n=6$ and $k=3$ from~\cite{KWZ}.}
    \label{fig:m=2}
\end{figure}

We now consider nonintersecting paths in networks $N_D$ for $D \in \Le(k,n)$ in order to understand bases for positroids.
An \emph{disjoint path system terminating at $J$} is a collection $(P_1, P_2, \dots, P_k)$ so that each $P_i$ is a path that starts at the $i$th source, there is no $i \neq j$ such that $P_i$ and $P_j$ share a vertex, and for each vertex of $u \in J$ there is one path that terminates at $u$.
Here, we allow $P_i$ to be the empty path that terminates at the $i$th source.
With this terminology, we have the following lemma.

\begin{lemma}[{\cite[Proposition 13]{Oh} and \cite[Theorem 6.5]{Pos}}]
A subset $J \subseteq [n]$ is a basis in the positroid associated with some $D \in \Le(k,n)$ if and only if there exists a disjoint path system terminating at $J$ in the network $N_D$.
\label{lem:Le}
\end{lemma}

We make a quick observation of the transversal matroid presented by $ \supp(\Meas(N))$ for any $N \in \Net(k,n)$ while noting this is not necessarily the same matroid as the linear matroid represented by $\Meas(N)$.

\begin{lemma}
If $N \in \Net(k,n)$, then $M[\supp(\Meas(N))]$ is a fundamental transversal matroid.
\label{lem:fundamental}
\end{lemma}
\begin{proof}
Take $N \in \Net(k,n)$ and let $\supp(\Meas(N)) = (S_1, S_2, \dots, S_k)$.
Let $a \in [n]$ be the $i$th source of $N$, then $a \in S_i$ and $a \not\in S_j$ for $j \neq i$.
It follows that $M[\supp(\Meas(N))]$ is a fundamental transversal matroid.
\end{proof}

\begin{figure}[h]
    \centering
    \begin{tikzpicture}
    \node[circle, fill=black, scale=0.5] (11) at (1,5) {};
    \node[circle, fill=black, scale=0.5] (21) at (1,4) {};
    \node[circle, fill=black, scale=0.5] (31) at (1,3) {};
    \node[circle, fill=black, scale=0.5] (41) at (1,2) {};
    \node[circle, fill=black, scale=0.5] (51) at (1,1) {};
    
    \node[circle, fill=black, scale=0.5] (12) at (2,5) {};
    \node[circle, fill=black, scale=0.5] (22) at (2,4) {};
    \node[circle, fill=black, scale=0.5] (32) at (2,3) {};
    \node[circle, fill=black, scale=0.5] (42) at (2,2) {};
    \node[circle, fill=black, scale=0.5] (52) at (2,1) {};

    \node[circle, fill=black, scale=0.5] (13) at (3,5) {};
    \node[circle, fill=black, scale=0.5] (23) at (3,4) {};
    \node[circle, fill=black, scale=0.5] (33) at (3,3) {};
    \node[circle, fill=black, scale=0.5] (43) at (3,2) {};
    \node[circle, fill=black, scale=0.5] (53) at (3,1) {};
    
    \node[circle, fill=black, scale=0.5] (14) at (4,5) {};
    \node[circle, fill=black, scale=0.5] (24) at (4,4) {};
    \node[circle, fill=black, scale=0.5] (34) at (4,3) {};
    \node[circle, fill=black, scale=0.5] (44) at (4,2) {};
    \node[circle, fill=black, scale=0.5] (54) at (4,1) {};
    
    \node[circle, fill=black, scale=0.5] (15) at (5,5) {};
    \node[circle, fill=black, scale=0.5] (25) at (5,4) {};
    \node[circle, fill=black, scale=0.5] (35) at (5,3) {};
    \node[circle, fill=black, scale=0.5] (45) at (5,2) {};
    \node[circle, fill=black, scale=0.5] (55) at (5,1) {};

    \draw[blue,   thick] (15) to (25);
    \draw[blue,   thick] (25) to[bend right=15] (24);
    \draw[blue,   thick] (24) to[bend right=15] (34);
    \draw[blue,   thick] (34) to[bend right=15] (33); 
    \draw[blue,   thick] (33) to[bend right=15] (43); 
    \draw[blue,   thick] (43) to[bend right=15] (42);
    \draw[blue,   thick] (42) to (41);
    
    \draw[red, densely dashdotted,  thick] (5.5,4) to (25);
    \draw[red, densely dashdotted,   thick] (25) to[bend left=15] (24);
    \draw[red, densely dashdotted,  thick] (24) to[bend left=15] (34);
    \draw[red, densely dashdotted,  thick] (34) to[bend left=15] (33); 
    \draw[red, densely dashdotted,  thick] (33) to[bend left=15] (43); 
    \draw[red, densely dashdotted,  thick] (43) to[bend left=15] (42);
    \draw[red, densely dashdotted,  thick] (42) to (52);

    \node (x) at (5.75,3) {};
    \node (y) at (7.25,3) {};
    \draw[-{latex}, ultra thick] (x)--(y);
    \node at (6.5,3.25) {flip};
    
    \begin{scope}[shift={(7,0)}]
       \node[circle, fill=black, scale=0.5] (11) at (1,5) {};
    \node[circle, fill=black, scale=0.5] (21) at (1,4) {};
    \node[circle, fill=black, scale=0.5] (31) at (1,3) {};
    \node[circle, fill=black, scale=0.5] (41) at (1,2) {};
    \node[circle, fill=black, scale=0.5] (51) at (1,1) {};
    
    \node[circle, fill=black, scale=0.5] (12) at (2,5) {};
    \node[circle, fill=black, scale=0.5] (22) at (2,4) {};
    \node[circle, fill=black, scale=0.5] (32) at (2,3) {};
    \node[circle, fill=black, scale=0.5] (42) at (2,2) {};
    \node[circle, fill=black, scale=0.5] (52) at (2,1) {};

    \node[circle, fill=black, scale=0.5] (13) at (3,5) {};
    \node[circle, fill=black, scale=0.5] (23) at (3,4) {};
    \node[circle, fill=black, scale=0.5] (33) at (3,3) {};
    \node[circle, fill=black, scale=0.5] (43) at (3,2) {};
    \node[circle, fill=black, scale=0.5] (53) at (3,1) {};
    
    \node[circle, fill=black, scale=0.5] (14) at (4,5) {};
    \node[circle, fill=black, scale=0.5] (24) at (4,4) {};
    \node[circle, fill=black, scale=0.5] (34) at (4,3) {};
    \node[circle, fill=black, scale=0.5] (44) at (4,2) {};
    \node[circle, fill=black, scale=0.5] (54) at (4,1) {};
    
    \node[circle, fill=black, scale=0.5] (15) at (5,5) {};
    \node[circle, fill=black, scale=0.5] (25) at (5,4) {};
    \node[circle, fill=black, scale=0.5] (35) at (5,3) {};
    \node[circle, fill=black, scale=0.5] (45) at (5,2) {};
    \node[circle, fill=black, scale=0.5] (55) at (5,1) {};

    \draw[blue,   thick] (15)--(14)--(13)--(12)--(11)--(21)--(31)--(41);  
    
    \draw[red, densely dashdotted,  thick] (5.5,4)--(25)--(24)--(34)--(33)--(43)--(42)--(52);
    \end{scope}
    
    \begin{scope}[shift={(0,-6)}]
    \node[circle, fill=black, scale=0.5] (11) at (1,5) {};
    \node[circle, fill=black, scale=0.5] (21) at (1,4) {};
    \node[circle, fill=black, scale=0.5] (31) at (1,3) {};
    \node[circle, fill=black, scale=0.5] (41) at (1,2) {};
    \node[circle, fill=black, scale=0.5] (51) at (1,1) {};
    
    \node[circle, fill=black, scale=0.5] (12) at (2,5) {};
    \node[circle, fill=black, scale=0.5] (22) at (2,4) {};
    \node[circle, fill=black, scale=0.5] (32) at (2,3) {};
    \node[circle, fill=black, scale=0.5] (42) at (2,2) {};
    \node[circle, fill=black, scale=0.5] (52) at (2,1) {};

    \node[circle, fill=black, scale=0.5] (13) at (3,5) {};
    \node[circle, fill=black, scale=0.5] (23) at (3,4) {};
    \node[circle, fill=black, scale=0.5] (33) at (3,3) {};
    \node[circle, fill=black, scale=0.5] (43) at (3,2) {};
    \node[circle, fill=black, scale=0.5] (53) at (3,1) {};
    
    \node[circle, fill=black, scale=0.5] (14) at (4,5) {};
    \node[circle, fill=black, scale=0.5] (24) at (4,4) {};
    \node[circle, fill=black, scale=0.5] (34) at (4,3) {};
    \node[circle, fill=black, scale=0.5] (44) at (4,2) {};
    \node[circle, fill=black, scale=0.5] (54) at (4,1) {};
    
    \node[circle, fill=black, scale=0.5] (15) at (5,5) {};
    \node[circle, fill=black, scale=0.5] (25) at (5,4) {};
    \node[circle, fill=black, scale=0.5] (35) at (5,3) {};
    \node[circle, fill=black, scale=0.5] (45) at (5,2) {};
    \node[circle, fill=black, scale=0.5] (55) at (5,1) {};

    \draw[blue,   thick] (15) to (25);
    \draw[blue,   thick] (25) to[bend right=15] (35);
    \draw[blue,   thick] (35) to[bend right=15] (34);
    \draw[blue,   thick] (34) to[bend right=15] (44); 
    \draw[blue,   thick] (44) to[bend right=15] (43); 
    \draw[blue,   thick] (43) to[bend right=15] (53);
    \draw[blue,   thick] (53) to (3,0.5);
    
    \draw[red, densely dashdotted,  thick] (5.5,4) to (25);
    \draw[red, densely dashdotted,  thick] (25) to[bend left=15] (35);
    \draw[red, densely dashdotted, thick]  (35) to[bend left=15] (34);
    \draw[red, densely dashdotted,  thick] (34) to[bend left=15] (44); 
    \draw[red, densely dashdotted,  thick] (44) to[bend left=15] (43); 
    \draw[red, densely dashdotted,  thick] (43) to[bend left=15] (53);
    \draw[red, densely dashdotted,  thick] (53) to (52);

    \node (x) at (5.75,3) {};
    \node (y) at (7.25,3) {};
    \draw[-{latex}, ultra thick] (x)--(y);
    \node at (6.5,3.25) {flip};
    
    \begin{scope}[shift={(7,0)}]
       \node[circle, fill=black, scale=0.5] (11) at (1,5) {};
    \node[circle, fill=black, scale=0.5] (21) at (1,4) {};
    \node[circle, fill=black, scale=0.5] (31) at (1,3) {};
    \node[circle, fill=black, scale=0.5] (41) at (1,2) {};
    \node[circle, fill=black, scale=0.5] (51) at (1,1) {};
    
    \node[circle, fill=black, scale=0.5] (12) at (2,5) {};
    \node[circle, fill=black, scale=0.5] (22) at (2,4) {};
    \node[circle, fill=black, scale=0.5] (32) at (2,3) {};
    \node[circle, fill=black, scale=0.5] (42) at (2,2) {};
    \node[circle, fill=black, scale=0.5] (52) at (2,1) {};

    \node[circle, fill=black, scale=0.5] (13) at (3,5) {};
    \node[circle, fill=black, scale=0.5] (23) at (3,4) {};
    \node[circle, fill=black, scale=0.5] (33) at (3,3) {};
    \node[circle, fill=black, scale=0.5] (43) at (3,2) {};
    \node[circle, fill=black, scale=0.5] (53) at (3,1) {};
    
    \node[circle, fill=black, scale=0.5] (14) at (4,5) {};
    \node[circle, fill=black, scale=0.5] (24) at (4,4) {};
    \node[circle, fill=black, scale=0.5] (34) at (4,3) {};
    \node[circle, fill=black, scale=0.5] (44) at (4,2) {};
    \node[circle, fill=black, scale=0.5] (54) at (4,1) {};
    
    \node[circle, fill=black, scale=0.5] (15) at (5,5) {};
    \node[circle, fill=black, scale=0.5] (25) at (5,4) {};
    \node[circle, fill=black, scale=0.5] (35) at (5,3) {};
    \node[circle, fill=black, scale=0.5] (45) at (5,2) {};
    \node[circle, fill=black, scale=0.5] (55) at (5,1) {};

    \draw[blue,  thick] (15)--(14)--(13);
    \draw[red, densely dashdotted, thick] (13)--(12);
    \draw[blue, thick] (12)--(22)--(32)--(42)--(52);  
    
    \draw[red, densely dashdotted,   thick] (5.5,4)--(25)--(35)--(34)--(44)--(43)--(53);
    \draw[blue, thick] (53)--(3,0.5);
    \end{scope}
    \end{scope}
    \end{tikzpicture}
    \caption{Flips used to remove an intersection.}
    \label{fig:flip}
\end{figure}

We now prove our main theorem, where we restrict the $\Sq$-diagrams.
Here we find that the transversal matroid presented by $\supp(\Meas(N_D))$ and the linear matroid represented by $\Meas(N_D)$ coincide when $D \in \Sq(k,n)$.

\begin{theorem}
If $D \in \Sq(k,n)$, then the positroid $M_{\Meas(N_D)}$ is a fundamental transversal matroid satisfying
\[M_{\Meas(N_D)} = M[\supp(\Meas(N_D))]\]
where $\supp(\Meas(N_D))$ is a noncrossing minimal presentation.
\label{thm:Sq}
\end{theorem}

\begin{proof}
We consider the network $N_D$ for some $D \in \Sq(k,n)$ and take $J \subseteq [n]$ with $|J| = k$.
By Lemma~\ref{lem:Le} $J$ is a basis of the positriod represented by $\Meas(N_D)$ if and only if there is a disjoint path system in $N_D$ terminating at $J$.
Also, $J$ is a base of the transversal matroid presented by $\supp(\Meas(N_D))$ if and only if there is a collection of paths $(P_1, P_2, \dots, P_k)$ terminating at $J$ (with no conditions on being disjoint).
So, we must show that if there exists a collection of paths $(P_1, P_2, \dots, P_k)$ with $P_i$ starting at the $i$th source in $N_D$  where each member of $J$ has a path terminating at it, then there exists a disjoint path system in $N_D$ terminating at $J$.

Let $(P_1, P_2, \dots, P_k)$ be such a collection of paths terminating at $J$.
If two paths intersect, they must first meet at a vertex with one path coming from above and the other coming from the right.
Then the two may be identical for some time, but they must eventually diverge because no two paths in the collection terminate in the same place.
Now pick two paths that intersect so that the vertex they first meet in is as far to the bottom right as possible in the $k \times (n-k)$ box our $\Sq$-diagram lives inside.
It is then possible to perform an operation on the two paths that ``flip'' portions of the paths to remove this intersection.
Two cases of this operation are shown in Figure~\ref{fig:flip}, where this figure only shows the relevant portion of the paths near the intersection vertex.
There can be vertices present in the $\Sq$-network that are not explicitly shown in the figure (i.e. there can be vertices between the changes of direction that are not important to us).
However, starting with the vertices on the paths, the existence of each other vertex shown in the figure is implied by the definition of a $\Sq$-diagram.

There are four cases depending on whether the solid blue path diverges from the red by taking a horizontal or vertical step.
Also, the paths may agree by either first going downward or first going left.
The unpictured cases are analogous to those shown in Figure~\ref{fig:flip}.
In any case, new edges arising---and hence any new intersections---from the operation are strictly to the top and left side of the intersection we started with.
So, we can find a disjoint path system terminating at $J$ by iterating this operation.

We have shown that $M_{\Meas(N_D)} = M[\supp(\Meas(N_D))]$.
By Lemma~\ref{lem:MeasNoncrossing} we know that $\supp(\Meas(N_D))$ is noncrossing.
Also, by Lemma~\ref{lem:fundamental} we find that $M_{\Meas(N_D)} = M[\supp(\Meas(N_D))]$ is a fundamental transversal matroid.
It remains to verify that we are dealing with a minimal presentation.

Let $r_i \in \RR^n$ for $1 \leq i \leq k$ denote the rows of the boundary measurement matrix $\Meas(N_D)$.
Showing that we have a minimal presentation is equivalent to showing 
\[\supp(r_i) = [n] \sm H_i\]
for each $ 1 \leq i \leq k$ where $H_i$ is a hyperplane of the positroid.
We choose any $1 \leq i \leq k$  and set $H_i := [n] \sm \supp(r_i)$; so, $H_i$ contains all the indices where $r_i$ is $0$.
It is clear $\rank(H_i) < k$ because $H_i$ is missing the $i$th source along with everything the $i$th source is connected to by a path.
Hence, no subset of $H_i$ of size $k$ has a disjoint path system terminating at it.
Since $H_i$ contains the $k-1$ sources that are not the $i$th source, we find that $\rank(H_i) = k-1$.
The last thing to prove is that $H_i$ is a flat, meaning that if we add any element to $H_i$, its rank will increase.
Let us add any element $x \in [n] \sm H_i$ to $H_i$.
We take $S \cup \{x\} \subseteq H_i \cup \{x\}$ where $S$ is the set consisting of the $k-1$ sources in $H_i$.
Then $\rank(S \cup \{x\}) = k$ since there is a disjoint path system terminating at $S \cup \{x\}$.

%We will do this by citing geometric results.
%The geometric spaces used will not be defined since they will not be used outside of the next paragraph, but the interested reader can see the cited work for further details.
%Since $N$ is a $\Sq$-network, a source $i$ is connected to a sink $j$ if and only if $D$ has a $\bullet$ in the row corresponding to the source $i$ in the position corresponding to the column of sink $j$.
%So, $|D|$ is to equal the number of nonzero entries in $\Meas(N_D)$ minus $k$ since we have an entry equal to $1$ for each of the $k$ sources.
%We know that $|D|$ is the dimension of the space parameterized by the $\Le$-network $N_D$ in the totally nonnegative Grassmannian~\cite[Theorem 6.5]{Pos}.
%By~\cite[Theorem 5.9]{KLS} we obtain that $|D| = \dim(\Pi_D)$.
%Let $\cS = (S_1, S_2, \dots, S_k) = \supp(\Meas(N_D))$, then by~\cite[Theorem 3.2.1]{MarcottThesis} we see that $|D| = \dim(L(\cS))$ because $M_{\Meas(N_D)} = M[\supp(\Meas(N_D))]$ is a positroid so $\dim(\Pi_D) = \dim(L(\cS))$. 
%Now we have that 
%\[ |D| = -k + \sum_{i=1}^k |S_i|\]
%from which it follows that $\supp(\Meas(N_D))$ is a minimal presentation by using~\cite[Theorem 3.1.2]{MarcottThesis}.
\end{proof}

\section{Transversal Positroids Beyond $\Sq$-diagrams Via Cyclic Flats}

In this section, we find further families of transversal positroids by understanding their cyclic flats in certain cases.
Before looking at these families, we record some general results that will be of use.
Cyclic flats can be ordered by containment, and the resulting poset is a lattice.
A criterion determining if a given matroid is transversal was given by Mason~\cite{Mason} and later refined by Ingleton~\cite{Ingleton} using antichains in the lattice of cyclic flats and inequalities of ranks.
Bonin--Kung--de Mier~\cite{BoninKungdeMier} found an analogous result for fundamental transversal matroids.
The classification says that a matroid is transversal if and only if for all nonempty antichains $\cF$ of cyclic flats,
\begin{equation}\label{eq:trans}
    \rank(\cap \cF)\leq \sum_{\cF'\subseteq \cF}(-1)^{|\cF'|+1}\rank(\cup \cF')
    \tag{$\dagger$}
\end{equation}  
    where 
    \begin{align*}
        \cap\cF=\bigcap_{F\in \cF} F && \text{and} && \cup\cF' = \bigcup_{F \in \cF'}F.
    \end{align*}
Moreover, the matroid is a fundamental transversal matroid if and only if Equation~(\ref{eq:trans}) holds with equality for all nonempty antichains of cyclic flats.

We now make note of a few results that will be useful.
Whenever $\cF$ consists of a single cyclic flat, both sides of Equation~(\ref{eq:trans}) are equal.
We use the term \emph{trivial cyclic flat} to mean either the minimal or maximal element of the lattice of cyclic flats.
Since we only need to consider cases where $\cF$ is an antichain, we can always restrict our attention to nontrivial cyclic flats.
The \emph{width} of a poset is the size of the largest antichain.
If the lattice of cyclic flats of a matroid has width at most $2$, then the matroid is (bi)transversal\footnote{\emph{Bitransversal} means that both the matroid and its dual are transversal.}~\cite[Theorem 5.6]{BoninDMier}.
If $H$ is a cyclic hyperplane of a transversal matroid, then $E \setminus H$ is a member of any presentation of the matroid~\cite[Corollary 2.6]{ExtPres}.
It follows that a transversal matroid of rank $k$ has at most $k$ cyclic hyperplanes~\cite[Corollary 2.7]{ExtPres}.
Finally, we provide an alternative characterization of when a presentation of a transversal matroid is minimal.
A set system is a minimal presentation of the transversal matroid that it presents if and only if each set in the set system is a distinct cocircuit (i.e., complement of a hyperplane) of the transversal matroid~\cite[Theorem 3.2 and 3.7]{DB72}.

\subsection{Transversal Sparse Paving Positroids}
A rank $k$ matroid on ground set $E$ is \emph{sparse paving} if for every $A \in \binom{E}{k}$ either $A$ is a basis or $A$ is a circuit-hyperplane.
So, in a sparse paving matroid, the lattice of cyclic flats consists only of the circuit-hyperplanes along with the minimal and maximal elements. 
The following lemma computes the right side of Equation~(\ref{eq:trans}) in the case of an antichain of cyclic hyperplanes.

\begin{lemma}
    If $\cF \neq \emptyset$ contains only cyclic hyperplanes of a rank $k$ matroid, then
    \[\sum_{\cF'\subseteq \cF}(-1)^{|\cF'|+1}\rank(\cup \cF') = k - |\cF|.\]
    \label{lem:RHS}
\end{lemma}
\begin{proof}
    We start by observing that
    \[\rank(\cup \cF') = \begin{cases} 0 &|\cF'| = 0 \\ k-1 & |\cF'| = 1 \\  k & |\cF'| > 1\end{cases}\]
    since any hyperplane has rank equal to $k-1$ by definition, and the union of any two or more hyperplanes will be of full rank.
    So, we compute
    \begin{align*}
        \sum_{\cF'\subseteq \cF}(-1)^{|\cF'|+1}\rank(\cup \cF') &= |\cF|(k-1) + \sum_{j = 2}^{|\cF|}(-1)^{j+1}k\binom{|\cF|}{j}\\
        &= |\cF|(k-1) - (-k + |\cF|k)\\
        &= k - |\cF|
    \end{align*}
    where we have used the binomial theorem.    
\end{proof}

In~\cite{MachacekNasr}, the authors precisely characterize when a positroid is a sparse paving matroid.
We now describe which $\Le$-diagrams correspond to sparse paving positroids\footnote{We use $\Le$-diagrams here to match other sections of this work, but as in~\cite{MachacekNasr} one can also work with the \emph{Grassmann necklaces} or \emph{decorated permutations} for sparse paving positroids.}.
The uniform matroid $U_{k,n}$ is a positroid coming from the $\Le$-diagram whose shape is a $k \times (n-k)$ rectangle that is filled with a $\bullet$ position.
We will denote this $\Le$-diagram with $D_{k,n}$.
Assume $2 \leq k \leq n-2$ since the cases with rank $k \in \{0,1,n-1,n\}$ are trivial.
We label $n$ cells along the boundary of the $k \times (n-k)$ rectangle as follows:
\begin{enumerate}
    \item[(i)] we label the cell in the lower right corner with a $1$.
    \item[(ii)] we label, from right to left, the cells in the top row with $2, 3, \dots n - k+1$.
    \item[(iii)] we label, from top to bottom, the cells in the left-most column with $n-k+1, n-k+2, \dots, n$.
\end{enumerate}

The labeled positions are the only places a $\bullet$ can be removed in $D_{k,n}$ while still having a $\Le$-diagram.
For $A\subseteq [n]$ we define $D_{k,n}(A) \in \Le(k,n)$ to be the $\Le$-diagram obtained from $D_{k,n}$ by removing the $\bullet$ in the cell labeled $i$ for each $i \in A$.
When $i=1$, we modify the diagram by removing the $\bullet$ and the cell in that position, thus changing the underlying Young diagram.
In our notation $D_{k,n} = D_{k,n}(\emptyset)$.
See Figure~\ref{fig:le_numbering} for an example of the cell numbering and Figure~\ref{fig:le_sparse_pave} for examples of the construction. We find that $D_{k,n}(A)$ is typically not a $\Sq$-diagram; hence, we will indeed find new examples of transversal positroids compared to Theorem~\ref{thm:Sq}.

\begin{figure}[h]
     \begin{center}

         \begin{ytableau}
          7& 6 & 5 & 4 & 3 & 2\\
 8& & & & &\\
 9& & & & & \\
 10&& & & & 1
 \end{ytableau}

     \end{center}

     \caption{The labeling of the cells when $k = 4$ and $n = 10$.}
     \label{fig:le_numbering}
 \end{figure}

 \begin{figure}[h]
\begin{minipage}{.45\textwidth}
    \begin{center}
        \begin{ytableau}
\bullet& \bullet  & \bullet & \bullet &  & \bullet\\
\bullet&\bullet &\bullet &\bullet &\bullet &\bullet\\
\bullet& \bullet &\bullet &\bullet &\bullet &\bullet \\
&\bullet&\bullet&\bullet &\bullet  
\end{ytableau}
    \end{center}
\end{minipage}
\begin{minipage}{.45\textwidth}
    \begin{center}
        \begin{ytableau}
\bullet& \bullet  & \bullet & \bullet &  & \\
\bullet&\bullet &\bullet &\bullet &\bullet &\bullet\\
 & \bullet &\bullet &\bullet &\bullet &\bullet \\
\bullet&\bullet&\bullet&\bullet &\bullet  & \bullet
\end{ytableau}
    \end{center}
\end{minipage}

    \caption{On the left, we have the $\Le$-diagram $D_{4,10}(A)$, where $A=\{1,3,10\}$.
    On the right, we have $D_{4,10}(A)$ where $A=\{2,3,9\}$. }
    \label{fig:le_sparse_pave}
\end{figure}

To state a characterization of sparse paving positroid, we call $A\subseteq[n]$ \emph{non-adjacent} if for every $i\in A$ we have $i-1\notin A$ and $i+1\notin A$, where again we work modulo $n$.
Also recall, for any $0 \leq k \leq n$ and $i \in [n]$, we have defined the cyclic interval $C_{k,n}^{(i)} \in \binom{[n]}{k}$ as $C_{k,n}^{(i)} = \{i, i+1, \dots i + k -1\}$ where we work modulo $n$.

\begin{theorem}[{\cite[Theorem 5.4]{MachacekNasr}}]
    For $2 \leq k \leq n-2$, the positroid associated with the $\Le$-diagram $D_{k,n}(A)$ is sparse paving if and only if $A$ is non-adjacent.
    Moreover, the set of circuit hyperplanes of such a sparse paving positroid is $\{C_{k,n}^{(i)} : i \in A\}$.
    \label{thm:sparsepaving}
\end{theorem}

\begin{lemma}
    If $A \subseteq [n]$ is non-adjacent and nonempty, then \[\left|\bigcap_{i\in A}C^{(i)}_{k,n}\right| \leq \max \{0, k + 2 - 2|A|\}\]
    for any $2 \leq k \leq n-2$.
    \label{lem:intersection}
\end{lemma}
\begin{proof}
    We consider the complement in $[n]$ of the intersection in the lemma and wish to show
    \[\left|\bigcup_{i \in A} C^{(i+k)}_{n-k,n}\right| \geq \min \{n, n - k - 2 + 2|A|\}\]
    when $A \neq \emptyset$ is non-adjacent.
    The union of cyclic intervals will be a disjoint union of cyclic intervals.
    We concentrate on a single cyclic interval within the union.
    Assume $C^{(i_1)}_{n-k,n} \cup C^{(i_2)}_{n-k,n} \cup \cdots \cup C^{(i_{\ell})}_{n-k,n}$ is a cyclic interval starting at $i_1$.
    We know $|C^{(i_1)}_{n-k,n}|=n-k$, and because of the non-adjacent condition $|C^{(i_1)}_{n-k,n} \cup C^{(i_2)}_{n-k,n}|\geq n-k+2$ unless $n-k+2 > n$ and $C^{(i_1)}_{n-k,n} \cup C^{(i_2)}_{n-k,n} = [n]$.
    Each time we union the next cyclic interval, either our cardinality grows by at least $2$, or we will end up with all of $[n]$.
    The lemma now follows because any of the disjoint cyclic intervals in the union will have size at least $2$, and the same argument above will apply.
\end{proof}

\begin{theorem}
    Let $A \subseteq [n]$ be non-adjacent.
    For any $2 \leq k \leq n-2$, $D_{k,n}(A)$ is the $\Le$-digram of a transversal sparse paving positroid if and only if $|A| \leq k$.
    \label{thm:sparse_trans}
\end{theorem}

\begin{proof}
    Theorem~\ref{thm:sparsepaving} says $D_{k,n}(A)$ is the $\Le$-diagram for a sparse paving positroids and tells us exactly what the circuit hyperplanes are.
    We know that the positroid associated with $D_{k,n}(A)$ cannot be transversal if $|A| > k$ because a rank $k$ transversal matroid can have at most $k$ cyclic hyperplanes.
    We now check Equation~(\ref{eq:trans}), which we recall here: 
\[    \rank(\cap \cF)\leq \sum_{\cF'\subseteq \cF}(-1)^{|\cF'|+1}\rank(\cup \cF').\]
By  Lemma~\ref{lem:RHS}, for any nonempty collection of circuit hyperplanes we have
\[\sum_{\cF'\subseteq \cF}(-1)^{|\cF'|+1}\rank(\cup \cF')=k-|\cF|.\]
    And now by Lemma~\ref{lem:intersection}, along with the fact that $\rank(F) = |F|$ whenever $|F|<k$ since the smallest circuits have size $k$ in a rank $k$ sparse paving matroid, we have 
    \[\rank(\cap \cF) \leq \max\{0,k+2-2|\cF|\}\]
    when $|\cF|>1$, and $\rank(\cap \cF)=k-1$ when $\cF$ is a single circuit-hyperplane.
    Thus, the inequality in~(\ref{eq:trans}) holds whenever $|A|\leq k$.
    Therefore, we have shown exactly when a sparse paving positroid is a transversal matroid.
\end{proof}

\subsection{Transversal Positroids of Rank $2$}\label{sec:rank2}

Throughout this section, we only consider rank $2$ positroids. So, we will be dealing with $\Le$-diagrams $D$ fitting inside $2 \times (n-2)$ rectangle.
Hyperplanes will be rank $1$ flats.
A rank $2$ transversal matroid can have at most $2$ cyclic hyperplanes, and any matroid whose lattice of cyclic flats has width $2$ or less will be transversal.
Thus, a rank $2$ matroid is transversal if and only if it has at most $2$ cyclic hyperplanes.

It will be helpful in the statements of some results to assume positroids are loopless.
We can make this assumption without loss of generality because loops can be added and removed without affecting whether a matroid is a positroid or transversal.
From the $\Le$-diagram perspective, this amounts to adding or removing columns of the diagram that do not contain any $\bullet$'s.
In terms of minimal presentations with set systems, loops are just elements of the ground set that are not contained in any set in the set system.
An example of loop removal via $\Le$-diagrams is shown in Figure~\ref{fig:loopless_le_diagram_example}.

\begin{figure}[h]
    \begin{center}
        \begin{ytableau}
 \bullet & \bullet &  \bullet &  &    &  & \bullet  & \bullet  &   & \bullet \\
 \bullet &\bullet & \bullet &\bullet & & \bullet &\bullet 
\end{ytableau}
\hspace{1.5cm}
\begin{ytableau}
 \bullet & \bullet &  \bullet &     &  & \bullet  & \bullet  & \bullet \\
 \bullet &\bullet & \bullet &\bullet & \bullet &\bullet 
\end{ytableau}
    \end{center}
    \caption{A $\Le$-diagram of a rank $2$ positroid (left) and the $\Le$-diagram of the rank $2$ positroid obtain after removing loops.}
    \label{fig:loopless_le_diagram_example}
\end{figure}

Take a $\Le$-diagram $D$ fitting inside a $2 \times (n-2)$ rectangle, denote columns of $\Meas(N_D)$ by $\bv_i$ for $1 \leq i \leq n$.
Then $\{\bv_i : 1 \leq i \leq n\}$ can be taken to be the ground set of our positroid.
Hyperplanes are then the maximal collections of parallel column vectors, and a hyperplane is cyclic if there are at least $2$ parallel column vectors.

The standard basis vectors of $\RR^2$, which we denote $\be_1$ and $\be_2$, are always in our ground set as columns indexed by the two sources of the network.
The standard basis vector $\be_1$ is part of a cyclic hyperplane exactly when there is a sink in the network for which there is a path from the $1$st source to this sink but not from the $2$nd source.
Similarly, we have that $\be_2$ is in a cyclic hyperplane if and only if there is a sink in the network for which there a path from the $2$nd source but not the $1$st.
\begin{lemma}
    Let $D$ be a $\Le$-diagram of a loopless rank $2$ positroid.
    The standard basis vector $\be_1$ is part of a cyclic hyperplane if and only if there is a column in $D$ such that there is a $\bullet$ in the top row but not the bottom row.
\    \label{lem:e1}
\end{lemma}
\begin{proof}
    When there is a $\bullet$ only in the top row of a given column, this is the only way that the sink indexing this column is reachable only from the $1$st source. Thus, this is the only way to produce a cyclic hyperplane containing $\be_1$.
\end{proof}

\begin{lemma}
    Let $D$ be a $\Le$-diagram of a loopless rank $2$ positroid.
    The standard basis vector $\be_2$ is part of a cyclic hyperplane if and only if there is a column in $D$ such that there is a $\bullet$ in the bottom row, but not in the top row, and it occurs to the right of any column with $\bullet$'s in both rows.
    \label{lem:e2}
\end{lemma}
\begin{proof}
    If there is a column with $\bullet$'s in both rows, then for any column to the left, there is a path from the $1$st source to any sink with a $\bullet$ in the second row.
    Sinks below columns with a $\bullet$ only in the bottom row that are to the right of any column with $\bullet$'s in both rows are the only ways a sink can be reachable from the $2$nd source.
\end{proof}

\begin{definition}    
 For $i\geq 1$, let $\L(i)$ be the $\Le$-diagram
    \begin{center}
        \begin{ytableau}
 \empty &  &    \\
 \bullet &\bullet & \bullet
\end{ytableau}
$\cdots$
        \begin{ytableau}
 \empty &      & \bullet \\
  \bullet  & \bullet &\bullet 
\end{ytableau}
    \end{center}
where $i$ is the number of columns with only a single $\bullet$. 
\end{definition}

When $D$ is a $\Le$-diagram, we say that $D'$ is a \emph{contiguous subdiagram} of $D$ if $D'$ can be realized by taking consecutive columns in $D$. 
We say that an $\L(i)$ as a contiguous subdiagram of $D$ is \emph{maximal} if we cannot extend the contiguous subdiagram to the left to make an $\L(i+1)$.
The next lemma shows that having $\L(i)$ as a contiguous subdiagram produces a cyclic hyperplane.
This demonstrates a reason it will be convenient to assume we have no loops.
Inserting loops, or columns of unfilled cells in the $\Le$-diagram, could break a copy of $\L(i)$ as a contiguous subdiagram.
However, it will not change the linear dependence among the corresponding columns.
In Figure~\ref{fig:loopless_le_diagram_example}, the $\Le$-diagram on the right contains a maximal $\L(2)$.
The $\Le$-diagram on the left does not contain any $\L(2)$ since it is broken by the presence of a loop.
This is one reason we have made the assumption that our positroids are loopless.

\begin{lemma}\label{lem:Li}
    Let $D$ be a $\Le$-diagram of a loopless rank $2$ positroid.
    If $D$ contains a contiguous subdiagram of the form
    \begin{center}
        \begin{ytableau}
 \empty &  &    \\
 \bullet &\bullet & \bullet
\end{ytableau}
$\cdots$
        \begin{ytableau}
 \empty &      & \bullet \\
  \bullet  & \bullet &\bullet 
\end{ytableau}
    \end{center}
    (i.e., an $\L(i)$ for some $i \geq 1$), then all the columns of $\Meas(N_D)$ indexed by the sinks in the contiguous subdiagram belong to a cyclic hyperplane.
\end{lemma}
\begin{proof}
    Any path from either source to any of the sinks corresponding to columns in the picture in the lemma must pass through the vertex coming from the $\bullet$ in the bottom right corner.
    Hence, all these columns will be scalar multiples of each other with scalars coming from the weights along edges in the bottom row. 
\end{proof}

\begin{lemma}\label{lem:lin_ind}
    If a $\Le$-diagram $D$ of a loopless rank $2$ positroid has the following configuration
    \begin{center}
        \begin{ytableau}
\bullet\\
 \bullet
\end{ytableau}
$\cdots$
\begin{ytableau}
\bullet \\
 \bullet 
\end{ytableau}
\end{center}
where other columns can be filled in any way, respecting the $\Le$-diagram condition, then columns of $\Meas(N_D)$ indexed by the sinks below the two columns shown are linearly independent. 
\end{lemma}
\begin{proof}
Let $\bv_i$ and $\bv_j$ be the two column vectors on the right and left, respectively.
If
\[\bv_i = \begin{bmatrix} w_1 \\ w_2\end{bmatrix}\]
then 
\[\bv_j = \begin{bmatrix}\beta w_1 + \alpha \\ \beta w_2 \end{bmatrix}\]
where $\alpha \neq 0$.
Here $\beta$ is a product of edge weights along the bottom row, which is multiplied onto any paths leaving the bottom right $\bullet$ to the left.
Since there is a path going straight left from the $1$st source and turning down at the upper left $\bullet$ shown, it follows $\alpha \neq 0$.
\end{proof}

\begin{definition}
    For a loopless $\Le$-diagram $D$ of rank $2$ let $\tau_D$ be an integer defined in the following way: 
    \begin{enumerate}
        \item Set $\tau_D$ to be the number of distinct contiguous subdiagrams $\L(i)$ that are maximal, considering all values of $i \geq 1$. 
        \item Ignoring any columns covered by (1),
        \begin{enumerate}
            \item if there is a column with $\bullet$ in the only top row, add 1 to $\tau_D$; and 
            \item if there is a column with $\bullet$ in the only bottom row that is to the right of any column with $\bullet$'s in both rows, add 1 to $\tau_D$.
        \end{enumerate}
    \end{enumerate}
\end{definition}

\begin{figure}
    \begin{ytableau}
        \bullet & \bullet &         &         &         & \bullet &         & \bullet & \bullet\\
        \empty  &         & \bullet & \bullet & \bullet & \bullet & \bullet
    \end{ytableau}
    \caption{The $\Le$-diagram of a rank $2$ positroid on $[11]$ that has cyclic hyperplanes $\{1,2,3,10,11\}$, $\{4,5\}$, and $\{6,7,8,9\}$ each of which is a cyclic interval.}
    \label{fig:cyclic_flats}
\end{figure}

\begin{theorem}\label{thm:positroid_transv_rank2}
    Let $D$ be a $\Le$-diagram of a loopless rank $2$ positroid. The positroid associated with $D$ is transversal if and only if $\tau_{D} \leq 2$. Moreover, when the positroid is transversal it has a noncrossing minimal presentation.
    \label{thm:rk2}
\end{theorem}
\begin{proof}
        We know that the positroid is transversal if and only if there are at most $2$ cyclic hyperplanes.
    By Lemma~\ref{lem:e1}, Lemma~\ref{lem:e2}, and Lemma~\ref{lem:Li}, we find that the number of cyclic hyperplanes is at least as big as $\tau_D$.
    Lemma~\ref{lem:lin_ind} tells us $\tau_D$ is exactly the number of cyclic hyperplanes.
    By the loopless hypothesis, contiguous subdiagrams $\L(i)$ completely capture cyclic hyperplanes that do not contain either standard basis vector.
    It only remains to check that we have a noncrossing minimal presentation when $\tau_D \leq 2$.

    To obtain a minimal presentation, we will argue that all cyclic hyperplanes are cyclic intervals.
    The reader my find Figure~\ref{fig:cyclic_flats} helpful.
    If it exists, the cyclic hyperplane containing $\be_1$ must be a cyclic interval.
    The $1$st source is $1 \in [n]$ and is in this hyperplane.
    By the condition defining $\Le$-diagrams and the loopless hypothesis, any column containing a $\bullet$ only in the top row must make up contiguous subdiagrams of the leftmost or rightmost columns.
    If it exists, the cyclic hyperplane containing $\be_2$ must be an interval.
    Columns as described in Lemma~\ref{lem:e2} must come directly after the $2$nd source.
    Finally, a cyclic hyperplane coming from an $\L(i)$ must an interval as $\L(i)$ is a contiguous subdiagram.

    When $\tau_D = 0$ we have the uniform matroid $U_{2,n}$ which has a noncrossing minimal presentation.
    For example, we can take the presentation in Example~\ref{ex:uniform}.
    When $\tau_D = 1$, we have a cyclic hyperplane $Z$ and choose another hyperplane $H = \{i\}$ for any $i \in [n]\sm Z$.
    It is the case that $[n]\sm Z \neq \emptyset$ because $Z \neq [n]$ as our positroid is rank $2$.
    We claim $([n] \sm Z, [n]\sm H)$ is a noncrossing minimal presenation.
    When $\tau_D=2$ we have cyclic hyperplanes $Z_1$ and $Z_2$.
    We claim that $([n]\sm Z_1, [n] \sm Z_2)$ is a noncrossing minimal presenation.

    In any case, it is routine to check that we have given a presentation of the positroid as a transversal matroid.
    It is a minimal presentation because the sets in the set systems are complements of hyperplanes.
    Finally, the presentation is noncrossing because it consists of sets that are cyclic intervals.
    \end{proof}

\begin{example}
In Figure~\ref{fig:loopless_le_diagram_example}, we find that $\tau_D=2$ for the diagram $D$ on the right.
Hence, the positroid given by either diagram in this figure is transversal.
The diagram $D'$ in Figure~\ref{fig:cyclic_flats} has $\tau_{D'} = 3$ and thus the associated positroid is not transversal.
\end{example}

\section{Counterexamples to Noncrossing Minimal Presentations}

In this section, we provide a counterexample to a conjecture of Marcott~\cite[Conjecture 3.4.5]{MarcottThesis}, which claims that there exists a noncrossing minimal presentation for any positroid that is also a transversal matroid.
This conjecture is natural, since Marcott~\cite[Theorem 3.4.4]{MarcottThesis} has proven the converse and established special cases of the conjecture~\cite[Propositions 3.4.7 and 3.4.9]{MarcottThesis}.
Also, we have found in Theorem~\ref{thm:Sq} many examples of transversal positroids that have noncrossing minimal presentations due to the planarity of the networks parameterizing them.
In Lemma~\ref{lem:MeasNoncrossing} it was shown that $\supp(\Meas(N))$ is always noncrossing
However, the support of the boundary measurement matrix need not be a presentation for the corresponding positroid even if it is transversal.

\begin{theorem}
    If $n > 6$ and $A = \{2i+1 : 0 \leq i < n\}$, then the sparse paving positroid corresponding to $D_{n+1,2n}(A)$ is transversal and does not have a noncrossing minimal presentation.
    \label{thm:counter}
\end{theorem}
\begin{proof}
By Theorem~\ref{thm:sparse_trans}, the positroid corresponding to $D_{n+1,2n}(A)$ is transversal and thus must have minimal presentations.
Any minimal presentation consists of $n+1$ sets, each of which is the complement of a hyperplane.
The $n$ sets that are complements of the circuit hyperplanes $C^{(2i+1)}_{n+1,2n}$ must be in every presentation.
To show there is no noncrossing minimal presentation, we will demonstrate that one of the aforementioned sets in the presentation will cross the complement of any other hyperplane.
All other hyperplanes are independent sets of size $n$.
Each cyclic interval of length $n$ is contained in one of the circuit hyperplanes.
Every hyperplane that is not a circuit-hyperplane is an $n$-element independent set, and it cannot be a cyclic interval.

Take such a set
$H = \{x_1 < x_2 < \cdots < x_n\}$
and let $g_i = x_{i} - x_{i-1} - 1$ for $1 < i \leq n$ and $g_1 = (x_1+2n) - x_n - 1$ denote the \emph{gaps} between elements of $H$.
Since $H$ is not a cyclic interval, at least two gaps must be strictly positive.
Also,
\[\sum_i g_i = n\]
since the sum of the gaps plus $|H|=n$ must be $2n$. Now, select $i_0$ so that $g_{i_0}$ is the smallest strictly positive gap; it follows that $g_{i_0} \leq \frac{n}{2}$.

Now we let $a = x_{i_0}$ and think of hyperplane as
\[H = \{x_{i_0} <_a x_{i_0+1} <_a \cdots <_a x_{i_0-1}\}\]
with the cyclic order $<_a$. As $A$ only contains odd numbers, only one of $x_{i_0}+1$ or $x_{i_0}+2$ will be in $A$, and thus exactly one of $C^{(x_{i_0} + 1)}_{n+1,2n}$ or $C^{(x_{i_0} + 2)}_{n+1,2n}$ will be a circuit hyperplane. To this end, we let $Z := C^{(x_{i_0} + \epsilon)}_{n+1,2n}$ where $\epsilon \in \{1,2\}$ is selected so that $Z$ is in fact a circuit-hyperplane.
Hence, we find that $a \in H$ and $a \not\in Z$.

Next, let $b$ be any element such that $b \in Z$ and $b \not\in H$.
Since, $|Z| = |H|+1$ such a $b$ exists.
We let $c = x_{i_0-1}$, which satisfies
\[x_{i_0-1} = x_{i_0} + (n-2) + \sum_{i\neq i_0} g_i + 1\geq x_{i_0} + \frac{3n}{2} - 1\]
while
\[x_{i_0} + \epsilon + n\]
is the last element of the cyclic interval $Z = C^{(x_{i_0} + \epsilon)}_{n+1,2n}$.
Using $\epsilon \leq 2$ and $n > 6$, we find that
\[c>x_{i_0}+\epsilon+n,\]
and thus $c \not\in Z$ but $c \in H$.
Moreover, we find that $a <_a b <_a c$.

Finally, let $d = x_{i_0}-1$ and we see that $d \not\in H$ because $g_{i_0} > 0$.
Therefore, this choice of $a <_a b <_a c <_a d$ demonstrates that $[2n] \sm Z$ crosses $[2n] \sm H$.
Since $H$ was chosen arbitrarily from hyperplanes that are not circuit, the result is proven.
\end{proof}

\bibliographystyle{alphaurl}
\bibliography{refs}

\end{document}